\documentclass[10pt,a4paper]{amsart}
\usepackage{amsfonts}
\usepackage{epsfig}
\usepackage{graphicx}
\usepackage{amsmath}
\usepackage{amssymb}
\usepackage{txfonts}

\newcounter{main}

\newtheorem{theorem}{Theorem}[section]

\newtheorem{lemma}[theorem]{Lemma}
\newtheorem{corollary}[theorem]{Corollary}

\newtheorem{definition}{Definition}[section]
\newtheorem{maintheorem}{Theorem}

\newtheorem{maincorollary}{Corollary}

\newcommand{\blanksquare}{\,\,\,$\sqcup\!\!\!\!\sqcap$}

\newcounter{example}
\newenvironment{example}%
{{\stepcounter{example}}{\flushleft {\bf Example \arabic{example}:}}}%
{\par}

\def \dim{{\rm dim}}

\title[Non-uniform hyperbolicity for infinite dimensional cocycles]{Non-uniform hyperbolicity for infinite dimensional cocycles}

\author[M. Bessa]{M\'{a}rio Bessa}
\address{Departamento de Matem\'atica da Universidade do Porto,
Rua do Campo Alegre, 687,
4169-007 Porto, Portugal \\ ESTGOH-Instituto Polit\'ecnico de Coimbra, Rua General Santos Costa, 3400-124 Oliveira do Hospital, Portugal}
\email{bessa@fc.up.pt}
\author[M. Carvalho]{Maria Carvalho}
\address{Departamento de Matem\'atica da Universidade do Porto,
Rua do Campo Alegre, 687,
4169-007 Porto, Portugal}
\email{mpcarval@fc.up.pt}

\begin{document}

\begin{abstract}
Let $\mathcal{H}$ be an infinite dimensional separable Hilbert space, $X$ a compact Hausdorff space and $f:X \rightarrow X$ a 
homeomorphism which preserves a Borel ergodic measure which is positive on non-empty open sets. We prove that the non-uniformly 
Anosov cocycles are $C^0$-dense in the family of partially hyperbolic $f,\mathcal{H}$ - skew products with non-trivial unstable bundles.
\end{abstract}

\maketitle

\bigskip

\emph{keywords:} Random operators; dominated splitting; multiplicative ergodic theorem; Lyapunov exponents.\\
\normalsize

\begin{section}{Introduction} Let $\mathcal{H}$ be an infinite dimensional separable Hilbert space endowed with the inner
product $\langle\,,\, \rangle$, $\mathcal{C}(\mathcal{H})$ the set of linear compact operators acting in $\mathcal{H}$ with
the uniform norm, $f:X\rightarrow{X}$ a homeomorphism of a compact Hausdorff space $X$ and $\muup$ an $f$-invariant Borel 
ergodic measure which is positive on non-empty open subsets.

Given a family $(A_{x})_{x \in X}$ of operators in $\mathcal{C}(\mathcal{H})$, the associated skew product over $f$ and $\mathcal{H}$ (sometimes called cocycle)
is given by
$$
\begin{array}{cccc}
F(A): & X\times{\mathcal{H}} & \longrightarrow & X\times{\mathcal{H}} \\
& (x,v) & \longmapsto & (f(x),A(x)\cdot v).
\end{array}
$$
and so, for each $x\in X$ and $n \in \mathbb{N}_{0}$, we have $F^n(A)(x,v)=(f^n(x),A^n(x)\cdot v)$, where $(A^n)_{n \in \mathbb{N}_{0}}$
is the random sequence of linear maps defined by $A^{0}(x)=id$ and
$$A^{n}(x)=A(f^{n-1}(x))\circ...\circ A(f(x))\circ A(x).$$

We are interested in the asymptotic properties of the sequence $(((A^{*})^{n}A^{n}(x))^{\frac{1}{2n}})_{n \in \mathbb{N}}$, for $\muup$
almost every point $x$, where $A^{*}$ denotes the dual operator of $A$. Under the integrability condition
\begin{equation}\label{IC}
\int_{X}\log^{+}\|A(x)\| \, d\muup(x)<\infty,
\end{equation}
where $\log^{+}(y)=\text{max}\,\{0,\log(y)\}$, Ruelle's theorem (\cite{Ru}) gives, for $\muup$-almost every point $x\in X$,
a complete set of Lyapunov exponents of the above limit of operators and their associated invariant eigenspaces. In the sequel,
${C_{I}^{0}(X,\mathcal{C}(\mathcal{H}))}$ will stand for the continuous compact cocycles on $\mathcal{H}$ over $X$ satisfying
(\ref{IC}).

Reaching non-zero Lyapunov exponents is the key to attain other main informations about the geometric properties of a dynamical system (\cite{Ru}) and
the core of most of the results nowadays (\cite{T}). In \cite{BC}, the authors established generic properties of Oseledets-Ruelle's decompositions and
corresponding sets of exponents. More precisely, they proved that there exists a $C^{0}$-residual subset
$\mathcal{R}$ of ${C_{I}^{0}(X,\mathcal{C}(\mathcal{H}))}$ such that, for $A\in{\mathcal{R}}$ and $\muup$-almost every $x\in{X}$,
either the limit $\underset{n\rightarrow{\infty}}{\text{lim}}({A(x)^{*}}^{n}A(x)^{n})^{\frac{1}{2n}}$ is the null operator
or the Oseledets-Ruelle's splitting of $A$ along the orbit of $x$ is dominated (see definition \ref{Oseledets-dominated}). Naturally,
a dominated splitting does not prevent the existence of zero Lyapunov exponents. Moreover, in general, it is hard to remove them by a
small perturbation -- from \cite{B}, we know, for instance, that zero exponents are generic for non-hyperbolic
continuous cocycles of $\text{SL}(2,\mathbb{R})$) -- but this may be due to the lack of either regularity or hyperbolicity. And, indeed,
in \cite{BB}, a Lebesgue-preserving diffeomorphism $f$, exhibiting a \emph{partially hyperbolic} splitting, was perturbed, in the $C^1$-topology, 
to reach a positive \textbf{sum} of all the Lyapunov exponents along the central direction.

To extend their strategy to the infinite dimensional setting, we had to start finding an adequate definition of partial hyperbolicity and a way to
perturb within this context so that, while keeping invariant the sum of all the center-unstable Lyapunov exponents (our analogue to conservativeness), 
the sum of the central ones increases and no longer vanishes. Afterwards, inspired by \cite{BFP} and using techniques from \cite{BV} and \cite{BC}, another $C^0$-perturbation was performed so that \textbf{each} Lyapunov exponent of the central direction became different from zero.

\end{section}

\begin{section}{Definitions and elementary properties}

This section formulates the main hypothesis on the cocycle $A$.

\medskip

\begin{subsection}{Lyapunov exponents}\label{OR}

The following result ensures a (unique, see \cite{M}) Lyapunov-spectral decomposition for the limit of a random product of compact cocycles under the condition (\ref{IC}).

\begin{theorem}\label{Ruelle}(\cite[Corollary 2.2]{Ru})
Let $f:X\rightarrow X$ be a homeomorphism and $\muup$ an $f$-invariant Borel probability. If the cocycle $A$ belongs to
$C^{0}_{I}(X,\mathcal{C}(\mathcal{H}))$, then, for $\muup$-a.e. $x\in{X}$, we have the following properties:
\begin{enumerate}
\item [(a)] The limit
$\underset{n\rightarrow{\infty}}{\text{lim}}({A(x)^{*}}^{n}A(x)^{n})^{\frac{1}{2n}}$
exists and is a compact operator $\mathcal{L}(x)$.

\item [(b)] Let $e^{\lambda_{1}(x)}>e^{\lambda_{2}(x)}>...$ be the nonzero
eigenvalues of $\mathcal{L}(x)$ and $U_{1}(x)$, $U_{2}(x)$, ... the
associated eigenspaces whose dimensions are denoted by $n_{i}(x)$.
The sequence of real functions $\lambda_{i}(x)$, called
\textbf{Lyapunov exponents} of $A$, where $1\leq i(x) \leq j(x)$ and
$j(x)\in \mathbb{N} \cup \{\infty\}$, verifies:
\end{enumerate}

\begin{enumerate}
\item [(b.1)] The functions $\lambda_{i}(x)$, $i(x)$, $j(x)$ and $n_{i}(x)$
are $f$-invariant and depend in a measurable way on $x$.

\item [(b.2)] If $V_{i}(x)$ is the orthogonal complement of
$U_{1}(x)\oplus{U_{2}(x)}\oplus{...}\oplus{U_{i-1}(x)}$, for
$i<j(x)+1$, and $V_{j(x)+1}(x)=Ker(\mathcal{L}(x))$, then:
\begin{enumerate}

\item[(i)] $\underset{n\rightarrow{\infty}}{\text{lim}}\frac{1}{n}\log\|A^{n}(x)u\|=\lambda_{i}(x) \text{ if }u\in{V_{i}(x)\setminus V_{i+1}(x)}\text{ and } i<j(x)+1$;

\item[(ii)] $\underset{n\rightarrow{\infty}}{\text{lim}}\frac{1}{n}\log\|A^{n}(x)u\|=-\infty\text{ if }u\in{V_{j(x)+1}(x)}$.

\end{enumerate}
\end{enumerate}
\end{theorem}

\medskip

As we are assuming that $\muup$ is ergodic, the maps $i(x)$, $j(x)$, $n_{i}(x)$ and $\lambda_{i}(x)$ are constant
$\muup$-almost everywhere and so, to evaluate them, it is enough to consider a $\muup$-generic point. In what follows,
we will denote by $\mathcal{O}(A)$ the full $\muup$ measure set of points
given by this theorem. Notice that, since $\muup$ is positive on non-empty open subsets,
$\mathcal{O}(A)$ is dense in $X$.

\bigskip

The infinite dimension of $\mathcal{H}$ brings additional trouble while dealing with Oseledets-Ruelle's decompositions
because, in what follows, we will need to deal with finite-dimensional Oseledets' spaces. Next lemma aims at overcoming this difficulty.

\bigskip

\begin{lemma}\label{finite dimension}\cite[Lemma 3.3]{BC}
Let $A$ be an integrable compact operator and $\lambda_{i}(x)$, $U_{i}(x)$ as in Ruelle's theorem. If $\lambda_{i}(x)\neq -\infty$,
then $U_{i}(x)$ has finite dimension.
\end{lemma}
\end{subsection}

\bigskip

\begin{subsection}{Dominated sums}

\begin{definition}\label{domination}
Given $f$ and $A$ as above, a positive integer $\ell$ and a number $\alpha \in ]0,1[$, we say that a direct sum
$E_{1}(x)\oplus E_{2}(x)$, defined on an $f$-invariant set $K$, is \textit{\textbf{$\ell, \alpha$-dominated}}, a property we will denote
by $E_{1} \succ_{\ell,\alpha} E_{2}$, if
\begin{enumerate}
\item[(I)] $A(E_{i}(x)) \subset E_{i} (f(x))$ for every $x \in K$.
\item[(II)] The dimension of $E_{i}(x)$ is constant, for $i=1,2$.
\item[(III)] There is $\theta > 0$ such that, for every $x \in K$ and any pair of unit vectors
$u\in E_{2}(x)$ and $v\in E_{1}(x)$, one has
\begin{itemize}
\item[(III.1)] $\|A(x)\cdot v\| \geq \theta.$
\item[(III.2)] $\|A^{\ell}(x)\cdot u\|\leq \alpha \, \|A^{\ell}(x)\cdot v\|.$
\end{itemize}
\end{enumerate}
\end{definition}

\bigskip

Condition (III.2) is a standard hypothesis of the classical concept of domination, whereas (III.1) is an essential demand in the infinite dimensional context in
order to guarantee that $Ker(A(x))\cap E_1(x)=\{\vec{0}\}$, for all $x \in K$. Among finite dimensional automorphisms, domination implies that the angle between any two
subbundles of a dominated splitting is uniformly bounded away from zero -- a useful property while proving, for instance, that a dominated splitting extends
continuously to the boundary of the set where it is defined. Due to the lack of compactness of $\mathcal{O}(A)$ and the fact that we are dealing with compact operators acting on an infinite dimensional space -- so $A(x)$ is not invertible and its norm may not be uniformly bounded away from zero -- we cannot expect
such a strong statement in our setting, unless we relate, as we have done in (III.1), domination with non-zero norms.

\begin{definition}\label{finest}
An $\ell, \alpha$-dominated direct sum $E_{1}\oplus E_{2}$ on $X$ is said to be \textit{\textbf{the finest}} if, after a non-trivial decomposition of any of these two subspaces, the sum is no longer dominated.
\end{definition}

Given a finite dimensional dominated sum $E_{1}\oplus E_{2}$ on $X$, a finest one exists and is unique. But the continuation of a finest dominated sum is not necessarily finest for the perturbed cocycle. However, the set of cocycles to whom this happens is a closed meager set (\cite{BFP}) -- that is, its complement is an open and dense set, we denote by $\mathcal{SF}_{E_{1},E_{2}}$.

\begin{definition}\label{dominated splitting}
Given a positive integer $\ell$ and a number $\alpha \in \,\,]0,1[$, a splitting $E_{1}(x)\oplus E_{2}(x) \oplus \cdots \oplus\ E_{k}(x)=\mathcal{H}$, defined on an $f$-invariant set $K$, is \textit{$\ell, \alpha$-dominated} if, for any $x \in K$ and every $1 \leq i<j \leq k$, we have $E_{i}(x) \succ_{\ell,\alpha} E_{j}(x)$.
\end{definition}

\medskip

The splitting we are interested in is the one corresponding to the Lyapunov subspaces given by Ruelle's theorem (so, under our assumptions,
conditions (I) and (II) are fulfilled \emph{a priori}). In further sections we will always address to this specific decomposition.

\medskip

\begin{definition}\label{Oseledets-dominated}
The Oseledets-Ruelle's decomposition in $\mathcal{O}(A)$ is said to be \textit{$\ell, \alpha$-dominated} if there are a positive integer $\ell$ and a number $\alpha \in \,\, ]0,1[$ such that, for each $x \in \mathcal{O}(A)$, we may rewrite this decomposition as a direct sum of $k$ subspaces, say $E_{1}(x)\oplus E_{2}(x) \oplus \cdots \oplus E_{k}(x)=\mathcal{H}$ so that:
\begin{enumerate}
\item For all $i \in \{1,2,\cdots,k\}$, the dimension of $E_{i}(x)$ is independent of $x \in \mathcal{O}(A)$.
\item For all $x \in \mathcal{O}(A)$ and every $i<j$, we have $E_{i}(x) \succ_{\ell, \alpha} E_{j}(x)$.
\end{enumerate}
\end{definition}

\medskip

Notice that, by definition, a dominated Oseledets-Ruelle's decomposition respects the order of the Lyapunov exponents: for instance, $E_{1}$ is associated with a finite number of the first (biggest) Lyapunov exponents.

The main properties of the domination in finite dimensions may now be conveyed to a dominated Oseledets-Ruelle's decomposition
$$x \in \mathcal{O}(A) \mapsto E_{1}(x)\oplus E_{2}(x) \oplus \cdots \oplus E_{k}(x)$$
of any infinite-dimension cocycle.

\begin{lemma}\label{A is invertible in E1}\cite[Lemma 3.4]{BC}{  Partial inverses}

Take $x$ in $\mathcal{O}(A)$ and consider the first $m$ Lyapunov exponents,
say $\lambda_{1}\ > \lambda_{2}> ... > \lambda_{m}$, and $E_{1}(x)=U_{1}(x)\oplus U_{2}(x)\oplus...\oplus U_{m}(x)$ the corresponding subspace. If
$\lambda_{m}> -\infty$, then the restriction of the operator $A(x): E_{1}(x) \rightarrow E_{1}(f(x))$ is invertible and $A^{-1}_{f(x)}$
is compact.
\end{lemma}

\begin{lemma}\label{Transversality}\cite[Lemma 3.6]{BC}{  Transversality}

There exists a constant $\gamma \in \,\, ]0,\frac{\pi}{2}]$ such that, for any $x \in \mathcal{O}(A)$ and any disjoint subsets $I$ and $J$ of $\{1,2,\cdots,k\}$,
the angle between $\underset{i \in I}{\bigoplus} E_{i}(x)$ and $\underset{j \in J}{\bigoplus} E_{j}(x)$ is bigger than $\gamma$.
\end{lemma}

\begin{lemma}\label{Extension to the closure}\cite[Proposition 3.5]{BC}{  Extension to the closure}

An $\ell, \alpha$-dominated Oseledets-Ruelle's decomposition may be extended
continuously to an $\ell, \alpha$-dominated splitting over the closure of $\mathcal{O}(A)$.
\end{lemma}

In the sequel, we will assume that the Oseledets-Ruelle's decomposition associated to $A$ and defined in $\mathcal{O}(A)$
has been so extended to $X$.

\begin{lemma}\label{Persistence}\cite[Appendix B]{BDV}{  Persistence}

A dominated splitting persists under $C^0$-small perturbations of the cocycle within $C^{0}_{I}(X,\mathcal{C}(\mathcal{H}))$.
\end{lemma}

That is, given an infinite cocycle $A$ with an $\ell, \alpha$-dominated Oseledets-Ruelle's decomposition and an $\epsilon \in \, \, ]0, \alpha[$, there is $\delta>0$ such that, if $\|B-A\|<\delta$, then the Oseledets-Ruelle's decomposition of $B$ is $\ell, (\alpha-\epsilon)$-dominated and its subbundles have the same dimensions of the initial dominated splitting.
\end{subsection}

\begin{subsection}{Partial hyperbolicity}

\begin{definition}
A cocycle $A \in C^{0}_{I}(X,\mathcal{C}(\mathcal{H}))$ with an extended Oseledets-Ruelle's decomposition is said to be \textbf{partially hyperbolic} if, for any $x \in X$, this splitting may be rewritten as a direct sum of three subspaces, say $E_x^{u}\oplus E_x^{c}\oplus E_x^{s}=\mathcal{H}$, and there are $\ell \in \mathbb{N}$, $\alpha \in \,\, ]0,1[$ and $\beta \in \,\, ]0,1[$ such that:
\begin{enumerate}
\item $E_x^{u}$ contains only Oseledets' subspaces associated to positive Lyapunov exponents.
\item $E_x^{c}$ contains the sum of all the Oseledets' subspaces associated to zero Lyapunov exponents.
\item $E_x^{s}$ contains only Oseledets' spaces corresponding to negative Lyapunov exponents and includes the sum of all those spaces determined by Lyapunov exponents equal to $-\infty$.
\item $E^{u} \succ_{\ell, \alpha} E^{c}$ and $E^{c} \succ_{\ell, \beta} E^{s}$.
\end{enumerate}
\end{definition}

\medskip

Thus, the dynamics along $E^u$ (respectively $E^s$) is strongly expanding (resp. contracting), while the weaker forms of expansion or contraction are gathered inside $E^c$. Notice that, from Lemma~\ref{finite dimension}, we deduce that if $A$ is a partially hyperbolic cocycle,
then $E^{u}$ and $E^{cu}:=E^u\oplus E^c$ are finite-dimensional.

\begin{definition}
If, besides being partially hyperbolic, no Lyapunov exponent of $A$ is zero, we say that the cocycle is \textbf{non-uniformly Anosov}\footnote{This nomenclature was, to our knowledge, first used in cf.~\cite{An}}.
\end{definition}

\bigskip

\begin{example}\label{example}
Let $A\in C^{0}_{I}(X,\mathcal{C}(\mathcal{H}))$ be the constant cocycle given by
$$A(e_1,e_2,...,e_n,...)=\left(2e_1,e_2,\frac{e_3}{3},...,\frac{e_n}{n},...\right),$$
where $\mathcal{B}=\{e_n\}_{n\in\mathbb{N}}$ is a orthonormal basis of $\mathcal{H}$. The eigenspaces are one-dimensional and generated by
the elements of $\mathcal{B}$ with eigenvalues $2$, $1$ and, for $n\geq 3$, $\sigma_n=\frac{1}{n}$, respectively. Thus, for each $n\in\mathbb{N}$,
the Lyapunov exponents are equal to $\log(2)$, $0$ and, for $n\geq 3$, $\lambda_n=-\log(n)$. This cocycle is $\ell, \alpha$-partially hyperbolic with
$\ell=1$, $\alpha=\frac{1}{6}$, $\beta=\frac{1}{3}$, $E_x^u=\langle e_1 \rangle$, $E_x^c=\langle e_2 \rangle$ and $E_x^s=\langle \mathcal{B}\setminus \{e_1,e_2\} \rangle$.
\end{example}
\end{subsection}

\bigskip

From Lemma \ref{Persistence}, we conclude that

\begin{corollary}\label{Openness}
The set $\mathcal{PH}$ of partially hyperbolic cocycles is $C^0$-open in $C^{0}_{I}(X,\mathcal{C}(\mathcal{H}))$.
\end{corollary}

\medskip

\begin{subsection}{Center-unstable metric entropy}

It is known (\cite{O}) that, for $\muup$-almost everywhere $x$ in $\mathcal{O}(A)$, we have
\begin{equation}\label{det unstable}
\underset{n\rightarrow \infty}{\lim}\frac{1}{n}\log|\det(A^n(x))|_{E_x^{u}}|=\sum_{i=1}^{d}\lambda^i n_i,
\end{equation}
and
\begin{equation}\label{det center-unstable}
\underset{n\rightarrow \infty}{\lim}\frac{1}{n}\log|\det(A^n(x))|_{E_x^{cu}}|=\sum_{i=1}^{D}\lambda^i n_i,
\end{equation}
where $n_i$ is the dimension of $U^i(x)$, $\lambda^u$ and $\lambda^c$ represent Lyapunov exponents in
$E_x^{u}$ and $E_x^{c}$, and $d$ and $D$ are the dimensions of $E_x^{u}$ and $E_x^{cu}$, respectively. From now on, to simplify
the notation, $\mathcal{O}(A)$ also stands for this full $\muup$ measure set in $\mathcal{O}(A)$.

\bigskip

If the map $x \in \mathcal{O}(A) \rightarrow \log|\det(A(x))|_{E^{cu}}|$ is $\muup$-integrable, then, by Birkhoff's ergodic theorem, the latter sum is equal to
\begin{eqnarray*}
\sum_{i=1}^D\lambda^i&=&\underset{n\rightarrow \infty}{\lim}\frac{1}{n}\log|\det(A^n(x))|_{E_x^{cu}}|\\
&=&\underset{n\rightarrow \infty}{\lim}\frac{1}{n}\log\left[\,|\det(A(f^{n-1}(x))|_{E_{f^{n-1}(x)}^{cu}}|\times \cdots \times |\det(A(x)|_{E_x^{cu}}|\,\right]\\
&=&\underset{n\rightarrow \infty}{\lim}\frac{1}{n}\sum_{i=0}^{n-1}\log|\det(A(f^{i}(x))|_{E_{f^{i}(x)}^{cu}}|= \int_X\log|\det(A(x)|_{E_x^{cu}}|d\muup(x),
\end{eqnarray*}
which, in case the bundle $E^{cu}$ is associated only with non-negative Lyapunov exponents, is the classical formula for the metric entropy. This motivates the following concept.

\begin{definition}
The $\muup$ -- \textbf{center-unstable metric entropy} of a partially hyperbolic cocycle $A$ over $f$ is given by
\begin{equation}\label{entropy}
\hbar_{\muup}(A):=\sum_{i=1}^{D}\lambda^i n_i.
\end{equation}
\end{definition}

\end{subsection}

\end{section}

\medskip

\begin{section}{Statement of the results}

We say that a partially hyperbolic cocycle in $C^{0}_{I}(X,\mathcal{C}(\mathcal{H}))$ has \textbf{\emph{non-trivial unstable bundle}} if, for any $x \in X$,
the unstable space $E_x^u$ does not reduce to $\{\vec{0}\}$.

\begin{maintheorem}\label{Teorema1}
Let $A \in C^{0}_{I}(X,\mathcal{C}(\mathcal{H}))$ be a partially hyperbolic cocycle with a non-trivial unstable bundle and such that $\sum  n_{c,A}\, \lambda^c_A = 0$. Then, for any $\epsilon>0$, there exists a cocycle $B \in C^{0}_{I}(X,\mathcal{C}(\mathcal{H}))$ verifying:
\begin{enumerate}
\item[(a)] $B$ is partially hyperbolic.
\item[(b)] $B$ is $\epsilon$ - $C^0$-close to $A$.
\item[(c)] $\hbar_{\muup}(B)=\hbar_{\muup}(A)$.
\item[(d)] $\sum n_{c,B} \, \lambda^c_B >0$.
\end{enumerate}
\end{maintheorem}

\medskip

The perturbation needed for this result will diminish the strength of the expansion along $E^u$ while keeping the center-unstable entropy invariant, so the sum of all the central Lyapunov exponents will have to increase. However, if $D-d > 1$, this sum may still have terms equal to zero. To remove them, we need another perturbation that makes all these summands essentially equal without destroying the positiveness of the sum.

\medskip

\begin{maintheorem}\label{Teorema2}
Let $B\in C^{0}_{I}(X,\mathcal{C}(\mathcal{H}))$ be a partially hyperbolic cocycle with a non-trivial unstable bundle and such that $\sum n_{c,B} \, \lambda^c_B \neq 0$. Then, for any $\epsilon>0$, there exists a non-uniformly Anosov cocycle $C\in C^{0}_{I}(X,\mathcal{C}(\mathcal{H}))$ which is $\epsilon$- $C^0$-close to $B$.
\end{maintheorem}

\medskip

\begin{maincorollary}
Non-uniformly Anosov cocycles are $C^0$-dense in the subset of partially hyperbolic ones with non-trivial unstable bundles.
\end{maincorollary}

\medskip

\begin{example}
Consider again the cocycle $A$ of Example~\ref{example}. Given $\epsilon \in \, ]0,\frac{1}{2}[$, take the cocycle $C$ defined by
$$C(e_1,e_2,...,e_n,...)=\left(2e^{-\epsilon} e_1, e^\epsilon e_2,\frac{e_3}{3},...,\frac{e_n}{n},...\right).$$ Then $C$ is
partially hyperbolic with respect to the splitting of $A$, $\|C-A\| < \epsilon$, $\hbar_{\muup}(C)=[\log(2)-\epsilon)]+(\epsilon)]=\hbar_{\muup}(A)$
and $0 = \lambda^c_{A} < \lambda^c_{C} =\epsilon.$
\end{example}
\end{section}

\medskip

\begin{section}{Sketch of the proof}

Let $A\in C^{0}_{I}(X,\mathcal{C}(\mathcal{H}))$ be a partially hyperbolic cocycle and an extended Oseledets-Ruelle's decomposition $E^{u}\oplus E^{c}\oplus E^{s}=\mathcal{H}$, such that, for any $x \in \mathcal{O}(A)$, the space $E_{x,A}^{u}$ is non-trivial with dimension $d$ and the dimension of $E_x^{cu}$ is $D$. Then:
\begin{eqnarray*}
\sum n_{c,A} \, \lambda^c_A \neq 0 \,\,\, &\Rsh& \text{  proceed to section \ref{PC}}. \\
\sum n_{c,A} \, \lambda^c_A =0 \,\,\, &\Rsh& \text{  go to section \ref{PCU} and, with the result, head for section \ref{PC}}.
\end{eqnarray*}
\end{section}

\medskip

\begin{section}{Perturbation lemmas}\label{PL}

\begin{subsection}{Perturbation on the center-unstable space}\label{PCU}

We will start showing how to perturb a partially hyperbolic cocycle in order to increase the sum of the Lyapunov exponents along the central
directions without changing the center-unstable entropy invariant.

\medskip

Let $B(x,r)=\{y \in X \colon \text{dist}(x,y)< r\}$ denote the ball centered at $x$ with radius $r>0$. Recall that, as $\muup$ is ergodic and positive on nonempty open sets, there is a residual subset of $X$ whose elements have dense orbits by $f$.

\medskip

\begin{lemma}\label{Perturbation on the center-unstable space}
Let $A\in C_{I}^{0}(X,\mathcal{C}(\mathcal{H}))$ be a partially hyperbolic cocycle with $D=\dim(E^{cu})$ and $\sum n_{c,A} \, \lambda^c_A =0$. Fix a point $p \in \mathcal{O}(A)$ with dense orbit, an $r>0$ and an $\epsilon>0$ small enough. Then there exist $\delta>0$ and a cocycle $B\in C_{I}^{0}(X,\mathcal{C}(\mathcal{H}))$ such that $B$ is positive and
\begin{enumerate}
\item[(i)] $B(x)=A(x)$ if $x\in X \setminus B(p,r)$.
\item[(2i)] $B(x)\cdot v^s_x=A(x)\cdot v^s_x$, $\forall\, v^s_x\in E_{x,A}^{s}$ and $\forall\, x\in X$.
\item[(3i)] $B(x)\cdot v_x=A(x)\cdot R_x\cdot v_x$, $\forall v_x\in E_{x,A}^{cu}$ and $\forall x \in B(p,r)$, where $R_x$ belongs to the space $\text{SO}(D,E_{x,A}^{cu})$ of rotations in $\text{SL}(D,\mathbb{R})$.
\item[(4i)] $\|A-B\|\leq{\epsilon}$.
\item[(5i)] $\hbar_{\muup}(B)=\hbar_{\muup}(A)$.
\item[(6i)] $\sum n_{c,B} \, \lambda^c_{B} > 0.$
\end{enumerate}
\end{lemma}

\begin{proof}
Let $\etaup_p\colon X\rightarrow\mathbb{R}$ be a continuous map such that $|\etaup_p(\cdot)|\leq 1$, $\etaup_p(x)=0$ if $\text{dist}(x,p)\geq r$
and $\etaup_p(x)=1$ if $\text{dist}(x,p)\leq \varepsilon r$, for a chosen $\varepsilon \in ]0,1[$.

\medskip

Take $\delta:=\epsilon \, \|A_{E^{cu}}\|^{-1}$ (see Lemma~\ref{A is invertible in E1}) and consider an isotopy $\Phi \colon [0,1]\rightarrow \text{SO}(D,\mathbb{R})$
such that $0 < \|id-\Phi\| \leq \delta$, $\Phi(0)=id$ and $\Phi(1)$ is a rotation centered at $p$, $\delta$-$C^0$-close
to the $id$ but different from it.
\medskip

Given $x \in X$, we may write each vector $v_x$ of $\mathcal{H}_x=\mathcal{H}$ in a unique way as $v_x=v^{s}_x + v^{cu}_x$, where $v^{s}_x\in E_{x,A}^{s}$ and
$v^{cu}_x \in E_{x,A}^{cu}$. Fixing now $r>0$ and $B(p,r)$, define the cocycle $B$ by
$$B(x)\cdot (v_x)= A(x)\cdot (v_x^s) + A(x) \circ \Phi[\etaup_p(\text{dist}(p,x))]\cdot (v_x^{cu}).$$

Notice that, by construction, the spaces $E_{x,A}^{s}$ and $E_{x,A}^{cu}$ are $B$-invariant for all $x\in X$. And, as, for all $x\in X$, the map
$\Phi[\etaup_p(\text{dist}(p,x))]$ is a mere rotation acting in $E_{x,A}^{cu}$, we have $E_{x,A}^{s}=E_{x,B}^{s}$ and
$E_{x,A}^{cu}=E_{x,B}^{cu}$. However, for any point $x$ whose orbit visits $B(p,r)$, the dynamics along $E_{x,A}^{c}$ and $E_{x,A}^{u}$
changes after the perturbation.

\medskip

\noindent \textbf{Items (i), (2i) and (3i)}:

\medskip

Clearly the cocycle $B$ satisfies these properties.

\bigskip

\noindent \textbf{Item (4i)}:

\medskip

This follows directly from the choice of $\delta$, since
\begin{eqnarray*}
\|A-B\| &=& \underset{x \in X}{\sup} \,\|A(x)-B(x)\|=\underset{x \in X}{\sup} \, \underset{v \in \mathcal{H} \colon \|v\|=1}{\sup} \, \|A(x).v-B(x).v\|\\
&=&\underset{x \in X}{\sup} \, \,\underset{v \in \mathcal{H} \colon \|v\|=1}{\sup} \, \|A(x).v_x^{cu}-A(x) \circ \Phi[\etaup_p(\text{dist}(p,x))]\cdot v_x^{cu}\| \\
&\leq& \|A\| \, \, \|id-\Phi\| \leq \|A\| \, \,  \delta = \epsilon.
\end{eqnarray*}

\noindent Therefore, if $\epsilon$ is small enough, then $B$ has positive norm. Moreover,
$$\int_{X}\log^{+}\|B(y)\| \, d\muup(y)\leq \int_{X}\log^{+}\|A(y)\| \, d\muup(y) \, + \, \epsilon < \infty.$$
And so, in particular, $B$ has an Oseledets-Ruelle's decomposition on a set $\mathcal{O}(B)$. The perturbation used ensures that this splitting is
made of subspaces $E_B^s=E_A^s$, and $E_B^c$, $E_B^u$ whose dimensions, if $\delta$ is small enough, are equal to the corresponding ones of $A$ because these two finite dimensional bundles vary continuously with the cocycle (as asserted in section 4 of \cite{Ru}). So, on the full $\muup$ measure set $\mathcal{O}(A)\cap \mathcal{O}(B)$,
$E^s_B\oplus E^{cu}_B=E^s_A\oplus E^{cu}_A$. Moreover,

\medskip

\begin{corollary}
The splitting $\mathcal{H}=E^s_B\oplus E^{cu}_B=E^s_A\oplus E^{cu}_A$ is $\ell, \rho$-dominated for $B$, with $\rho=\frac{1+\alpha}{2}$.
\end{corollary}

\begin{proof}
Denote by $E_1$ and $E_2$ the spaces $E^{cu}_B=E^{cu}_A$ and $E^{s}_B=E^{s}_A$, respectively. By construction of $B$, conditions
(I) and (II) of Definition \ref{domination} are satisfied in $\mathcal{O}(B)$.

\medskip

By assumption, there is $\theta > 0$ such that, for every $x \in \mathcal{O}(A)$ and any unit vector $v_1\in E_{1}(x)$, one
has $\|A(x)\cdot v_1\| \geq \theta$. Then, if $\epsilon < \frac{\theta}{2}$, we have
$$\|B(x)\cdot v_1\| = \|(B-A)(x)\cdot v_1 + A(x) \cdot v_1\| \geq  \|A(x) \cdot v_1\| - \|(B-A)(x)\cdot v_1\| \geq \theta - \epsilon \geq \frac{1}{2}\theta.$$

\medskip

As $\|A-B\|< \epsilon$, we know that $\|B\| \leq \|A\|+ \epsilon$ and that, for any $m \in \mathbb{N}$, there is a constant $\mathcal{K}_m > 0$, depending only on $m$, $\|A\|$ and $\epsilon$, such that
$$\|A^m-B^m\|< \epsilon \, \mathcal{K}_m.$$
Therefore, given $x \in \mathcal{O}(A)$ and unit vectors $v_1\in E_{1}(x)$ and $v_2\in E_{2}(x)$, if
$$\epsilon \leq \frac{1}{2 \, \mathcal{K}_\ell}\frac{1-\alpha}{1+\alpha}\,\left(\frac{\theta}{2}\right)^\ell,$$
then
\begin{eqnarray*}
\|B^\ell(x)\cdot v_2\| &=&\|(B^\ell-A^\ell)(x)\cdot v_2 + A^\ell(x) \cdot v_2\| \\
&\leq& \|(B^\ell-A^\ell)(x)\cdot v_2\| + \|A^\ell(x) \cdot v_2\| \\
&\leq& \epsilon \, \mathcal{K}_\ell + \alpha \, \|A^\ell(x) \cdot v_1\| \\
&\leq& \epsilon \, \mathcal{K}_\ell + \alpha \, \|(A^\ell-B^\ell)(x)\cdot v_1\| + \alpha \, \|B^\ell(x) \cdot v_1\| \\
&\leq& (1+\alpha) \mathcal{K}_\ell \,\epsilon + \alpha \, \|B^\ell(x) \cdot v_1\| \\
&\leq& \frac{1-\alpha}{2}\left(\frac{\theta}{2}\right)^\ell + \alpha \, \|B^\ell(x) \cdot v_1\| \\
&\leq& \frac{1-\alpha}{2} \, \|B^\ell(x) \cdot v_1\| + \alpha \, \|B^\ell(x) \cdot v_1\| \\
&=& \rho \, \|B^\ell(x) \cdot v_1\|.
\end{eqnarray*}
\end{proof}
\medskip

\noindent \textbf{Item (5i)}:

\medskip

As $\Phi[\etaup_p(\text{dist}(p,x))]$ belongs to $\text{SO}(D, E_{x,A}^{cu})$ and $E_{x,B}^{cu}=E_{x,A}^{cu}$,
if $z$ belongs to $\mathcal{O}(A)\cap \mathcal{O}(B)$, we have
\begin{eqnarray*}
\hbar_{\muup}(B)&=&\underset{n\rightarrow \infty}{\lim}\frac{1}{n}\log|\det B^n(z)|_{E_{z,B}^{cu}}|= \underset{n\rightarrow \infty}{\lim}\frac{1}{n}\log|\det B(f^{n-1}(z))_{E_{f^{n-1}(z),B}^{cu}}|\times \cdots \times |\det B(z)_{E_{z,B}^{cu}}| \\
&=&\underset{n\rightarrow \infty}{\lim}\frac{1}{n}\log|\det A(f^{n-1}(z))_{E_{f^{n-1}(z),A}^{cu}}|\times \cdots \times |\det A(z)_{E_{z,A}^{cu}}| \\
&=&\underset{n\rightarrow \infty}{\lim}\frac{1}{n}\log|\det A^n(z)|_{E_{z,A}^{cu}}|=\hbar_{\muup}(A).
\end{eqnarray*}

\bigskip

\noindent \textbf{Item (6i)}:

\medskip

Finally, we have to estimate the sum $\sum n_{c,B} \lambda^c_{B}$ and prove that it is strictly positive. We will check that, for some $0 < \Delta < 1$, we have
$\sum n_{u,B} \, \lambda^u_{B} \simeq \sum n_{u,A} \lambda^u_{A} + d \ln(\Delta)$ and then argue with the already established invariance of the center-unstable entropy.

Keeping in mind the argument in \cite{BB}, notice that we do not perturb the base dynamics $f$ and also that, in spite of the infinite dimension of $\mathcal{H}$, we will
just have to perturb in a finite dimensional subspace of it.

\medskip

\begin{subsubsection}{\textbf{First case: $d=1$}}

Under this hypothesis, as $\muup$ is ergodic, to evaluate the unstable Lyapunov exponent of $B$ we are only due to determine the average growth of $B^n(x)v^u$ for any $x \in \mathcal{O}(A)\cap\mathcal{O}(B)$ and a unit vector $v^u \in E^u_{x,A}$, as indicated in (b.2) of Theorem \ref{Ruelle}.

\medskip

As $p$ is not periodic, we may consider $r$ small enough so that $B(p,r)\cap f^{-1}(B(p,r))=\emptyset$ and $B(p,r)\cap f(B(p,r))=\emptyset$. Take $x$ in $B(p,r)\cap\mathcal{O}(A)\cap\mathcal{O}(B)\cap\mathcal{P}$, where $\mathcal{P}$ is the full $\muup$ measure subset of $B(p,r)$ given by Poincar\'e's recurrence theorem. Given a unit vector $v^u$ in $E^u_{x,A}$, we have, by definition,
$$B(x)\cdot v^u=A(x)\circ \pi_{E^u_{x,A}} \circ \Phi[\etaup_p(\text{dist}(p,x))] \cdot v^u + A(x)\circ \pi_{E^c_{x,A}} \circ \Phi[\etaup_p(\text{dist}(p,x))] \cdot v^u$$
where $\pi_{E^u_{x,A}}$ is the projection on $E^u_{x,A}$ parallel to the bundle $E^c_{x,A}$ (and analogous definition for $\pi_{E^c_{x,A}}$). As $f(x) \notin B(p,r)$, then
\begin{eqnarray*}
B^2(x)\cdot v^u &=& \\
&=&A(f(x))\circ \pi_{E^u_{f(x),A}} \circ \Phi[\etaup_p(\text{dist}(p,f(x)))] \cdot A(x)\circ \pi_{E^u_{x,A}} \circ \Phi[\etaup_p(\text{dist}(p,x))] \cdot v^u \\
&+&A(f(x))\circ \pi_{E^u_{f(x),A}} \circ \Phi[\etaup_p(\text{dist}(p,f(x)))] \cdot A(x)\circ \pi_{E^c_{x,A}} \circ \Phi[\etaup_p(\text{dist}(p,x))] \cdot v^u \\
&+&A(f(x))\circ \pi_{E^c_{f(x),A}} \circ \Phi[\etaup_p(\text{dist}(p,f(x)))] \cdot A(x)\circ \pi_{E^u_{x,A}} \circ \Phi[\etaup_p(\text{dist}(p,x))] \cdot v^u \\
&+&A(f(x))\circ \pi_{E^c_{f(x),A}} \circ \Phi[\etaup_p(\text{dist}(p,f(x)))] \cdot A(x)\circ \pi_{E^c_{x,A}} \circ \Phi[\etaup_p(\text{dist}(p,x))] \cdot v^u
\end{eqnarray*}

\medskip

\noindent which reduces to
$$B^2(x)\cdot v^u=A^2(x)\circ \pi_{E^u_{x,A}} \cdot \Phi[\etaup_p(\text{dist}(p,x))] \cdot v^u + A^2(x)\circ \pi_{E^c_{x,A}} \cdot \Phi[\etaup_p(\text{dist}(p,x))] \cdot v^u.$$
In general, while the orbit keeps out of $B(p,r)$, we have
$$B^j(x)\cdot v^u=A^j(x)\circ \pi_{E^u_{x,A}} \cdot \Phi[\etaup_p(\text{dist}(p,x))] \cdot v^u + A^j(x)\circ \pi_{E^c_{x,A}} \cdot \Phi[\etaup_p(\text{dist}(p,x))] \cdot v^u.$$
Now, by the domination of $E^u$ over $E^c$ under the action of $A$, there are constants $C$ and $\beta \in \, ]0,1[$ such that, for any $z \in X$, any $m \in \mathbb{N}$ and any unit vectors $w_1 \in E^u$ and $w_2 \in E^c$, we have
\begin{equation}\label{d}
\|A^m \cdot w_2\|\leq C\beta^m \|A^m \cdot w_1\|.
\end{equation}
Thus, the first component of $B^j(x)\cdot v^u$ dominates the second and contributes to $\lambda^u_B$ with an approximate rate of
\begin{equation}\label{lu}
\frac{1}{j}\log \|B^j(x)\cdot v^u\| \sim \lambda^u_A + log(\Delta) < \lambda^u_A,
\end{equation}
where $\Delta \in \, ]0,1[$ is an upper bound of the set $\{\cos(\omega_z): z \in X \}$ and $\omega_z$ stands for the small angle of rotation displaced by the action of $\Phi[\etaup_p(\text{dist}(p,z))]$.

\medskip

This component of the orbit out of $B(p,r)$ is the one that suits better our purposes to decrease the unstable Lyapunov exponent, and so it is essential that it lasts more than the visits to $B(p,r)$. For that, we demand that $r$ is small in order to guarantee, by Kac's theorem, that the expected first return to $B(p,r)$ of the orbit of $x$ amounts to the rather big fraction $\frac{1}{\muup \left(B(p,r)\right)}$, and so the contribution of the piece of orbit out of $B(p,r)$ for the estimation (\ref{lu}) is the prevalent one.

\medskip

Meanwhile, we know that, with $\muup$  probability one, the orbit of $x$ will eventually return to $B(p,r)$. If that happens first at an iterate $N$, then
\begin{eqnarray*}
\lefteqn{B^{N+1}(x)\cdot v^u=} \\
&=& B(f^N(x))\left(A^N(x) \circ \pi_{E^u_{x,A}} \cdot \Phi[\etaup_p(\text{dist}(p,x))] \cdot v^u\right)+ \\
& & B(f^N(x))\left(A^N(x) \circ \pi_{E^c_{x,A}} \cdot \Phi[\etaup_p(\text{dist}(p,x))] \cdot v^u\right)= \\
&=& A(f^N(x))\circ \pi_{E^u_{f^N(x),A}} \circ \Phi[\etaup_p(\text{dist}(p,f^N(x)))] \cdot \left(A^N(x)\circ \pi_{E^u_{x,A}} \circ \Phi[\etaup_p(\text{dist}(p,x))] \cdot v^u\right)+ \\
& & A(f^N(x))\circ \pi_{E^u_{f^N(x),A}} \circ \Phi[\etaup_p(\text{dist}(p,f^N(x)))] \cdot \left(A^N(x)\circ \pi_{E^c_{x,A}} \circ \Phi[\etaup_p(\text{dist}(p,x))] \cdot v^u\right)+ \\
& & A(f^N(x))\circ \pi_{E^c_{f^N(x),A}} \circ \Phi[\etaup_p(\text{dist}(p,f^N(x)))] \cdot \left(A^N(x)\circ \pi_{E^u_{x,A}} \circ \Phi[\etaup_p(\text{dist}(p,x))] \cdot v^u\right)+ \\
& & A(f^N(x))\circ \pi_{E^c_{f^N(x),A}} \circ \Phi[\etaup_p(\text{dist}(p,f^N(x)))] \cdot \left(A^N(x)\circ \pi_{E^c_{x,A}} \circ \Phi[\etaup_p(\text{dist}(p,x))] \cdot v^u\right).
\end{eqnarray*}

\bigskip
\bigskip

We proceed arguing as before, using domination, to see that the first of these four summands controls the others and yields an average as in (\ref{lu}).

\bigskip
\bigskip

\noindent (I) To compare the second term with the first, we apply the domination property (\ref{d}) to
$$w_1=\pi_{E^u_{f^N(x),A}} \circ \Phi[\etaup_p(\text{dist}(p,f^N(x)))] \cdot \left(A^N(x)\circ \pi_{E^u_{x,A}} \circ \Phi[\etaup_p(\text{dist}(p,x))] \cdot v^u\right) \in E^u$$
\noindent and
$$w_2=\pi_{E^c_{f^N(x),A}} \circ \Phi[\etaup_p(\text{dist}(p,f^N(x)))] \cdot \left(A^N(x)\circ \pi_{E^c_{x,A}} \circ \Phi[\etaup_p(\text{dist}(p,x))] \cdot v^u\right) \in E^c$$

\noindent in order to get

$$\frac{\|A(f^N(x))\circ \pi_{E^c_{f^N(x),A}} \circ \Phi[\etaup_p(\text{dist}(p,f^N(x)))] \cdot \left(A^N(x)\circ \pi_{E^c_{x,A}} \circ \Phi[\etaup_p(\text{dist}(p,x))] \cdot v^u\right)\|}{\|A(f^N(x))\circ \pi_{E^u_{f^N(x),A}} \circ \Phi[\etaup_p(\text{dist}(p,f^N(x)))] \cdot \left(A^N(x)\circ \pi_{E^u_{x,A}} \circ \Phi[\etaup_p(\text{dist}(p,x))] \cdot v^u\right)\|} \leq C\beta.$$

\noindent This term is controlled, in the limit, taking the logarithm and dividing by the iterate.

\bigskip
\bigskip

\noindent (II) In a similar way, we compare
$$w_1=\pi_{E^u_{f^N(x),A}} \circ \Phi[\etaup_p(\text{dist}(p,f^N(x)))] \cdot \left(A^N(x)\circ \pi_{E^u_{x,A}} \circ \Phi[\etaup_p(\text{dist}(p,x))] \cdot v^u\right) \in E^u$$
and
$$w_2=\pi_{E^c_{f^N(x),A}} \circ \Phi[\etaup_p(\text{dist}(p,f^N(x)))] \cdot \left(A^N(x)\circ \pi_{E^u_{x,A}} \circ \Phi[\etaup_p(\text{dist}(p,x))] \cdot v^u\right) \in E^c$$
and obtain

$$\frac{\|A(f^N(x))\circ \pi_{E^c_{f^N(x),A}} \circ \Phi[\etaup_p(\text{dist}(p,f^N(x)))] \cdot \left(A^N(x)\circ \pi_{E^u_{x,A}} \circ \Phi[\etaup_p(\text{dist}(p,x))] \cdot v^u\right)\|}{\|A(f^N(x))\circ \pi_{E^u_{f^N(x),A}} \circ \Phi[\etaup_p(\text{dist}(p,f^N(x)))] \cdot \left(A^N(x)\circ \pi_{E^u_{x,A}} \circ \Phi[\etaup_p(\text{dist}(p,x))] \cdot v^u\right)\|} \leq C\beta.$$

\bigskip
\bigskip

\noindent (III) Concerning the quotient
\begin{equation}\label{III}
\frac{\|A(f^N(x))\circ \pi_{E^u_{f^N(x),A}} \circ \Phi[\etaup_p(\text{dist}(p,f^N(x)))] \cdot \left(A^N(x)\circ \pi_{E^c_{x,A}} \circ \Phi[\etaup_p(\text{dist}(p,x))] \cdot v^u\right)\|}{\|A(f^N(x))\circ \pi_{E^u_{f^N(x),A}} \circ \Phi[\etaup_p(\text{dist}(p,f^N(x)))] \cdot \left(A^N(x)\circ \pi_{E^u_{x,A}} \circ \Phi[\etaup_p(\text{dist}(p,x))] \cdot v^u\right)\|}
\end{equation}
\noindent notice that
$$w_1=\pi_{E^u_{x,A}} \circ \Phi[\etaup_p(\text{dist}(p,x))] \cdot v^u \in E^u$$
and
$$w_2=\pi_{E^c_{x,A}} \circ \Phi[\etaup_p(\text{dist}(p,x))] \cdot v^u \in E^c$$

\noindent so, again by the domination inequality (\ref{d}),

$$\|A^N(x)\circ \pi_{E^c_{x,A}} \circ \Phi[\etaup_p(\text{dist}(p,x))] \cdot v^u\| \leq C\beta^N \|A^N(x)\circ \pi_{E^u_{x,A}} \circ \Phi[\etaup_p(\text{dist}(p,x))] \cdot v^u\|.$$

\noindent Thus, as $\Phi[\etaup_p(\text{dist}(p,f^N(x)))]$ is an isometry,
$$\|\Phi[\etaup_p(\text{dist}(p,f^N(x)))] \left(A^N(x)\circ \pi_{E^c_{x,A}} \circ \Phi[\etaup_p(\text{dist}(p,x))] \cdot v^u\right)\| \leq $$
$$C\beta^N \|\Phi[\etaup_p(\text{dist}(p,f^N(x)))]\left(A^N(x)\circ \pi_{E^u_{x,A}} \circ \Phi[\etaup_p(\text{dist}(p,x))] \cdot v^u\right)\|.$$

\medskip

\noindent Applying now to both vectors the projection $\pi_{E^u_{f^N(x),A}}$, whose norm is uniformly bounded in $B(p,r)$ by a constant $C_1$ (see Lemma~\ref{Transversality}), we obtain

$$\|\pi_{E^u_{f^N(x),A}} \circ \Phi[\etaup_p(\text{dist}(p,f^N(x)))] \left(A^N(x)\circ \pi_{E^c_{x,A}} \circ \Phi[\etaup_p(\text{dist}(p,x))] \cdot v^u\right)\| \leq $$ $$C_1C\beta^N \|\pi_{E^u_{f^N(x),A}} \circ \Phi[\etaup_p(\text{dist}(p,f^N(x)))]\left(A^N(x)\circ \pi_{E^u_{x,A}} \circ \Phi[\etaup_p(\text{dist}(p,x))] \cdot v^u\right)\|.$$

\medskip

\noindent The resulting vectors, appearing on the two sides of this inequality, say
$$\vartheta_1=\pi_{E^u_{f^N(x),A}} \circ \Phi[\etaup_p(\text{dist}(p,f^N(x)))] \left(A^N(x)\circ \pi_{E^c_{x,A}} \circ \Phi[\etaup_p(\text{dist}(p,x))] \cdot v^u\right)$$
and
$$\vartheta_2=\pi_{E^u_{f^N(x),A}} \circ \Phi[\etaup_p(\text{dist}(p,f^N(x)))]\left(A^N(x)\circ \pi_{E^u_{x,A}} \circ \Phi[\etaup_p(\text{dist}(p,x))] \cdot v^u\right),$$

\noindent both live in the one-dimensional subspace $E^u_{f^N(x),A}$, so there is a scalar $a$ such that $\vartheta_2=a \cdot \vartheta_1$, with $|a|\leq C\beta^N$. Hence, applying $A(f^N(x))$, we reach the inequality

$$\|A(f^N(x)) \circ \pi_{E^u_{f^N(x),A}} \circ \Phi[\etaup_p(\text{dist}(p,f^N(x)))] \left(A^N(x)\circ \pi_{E^c_{x,A}} \circ \Phi[\etaup_p(\text{dist}(p,x))] \cdot v^u\right)\|=$$
$$\|a \cdot A(f^N(x)) \circ \pi_{E^u_{f^N(x),A}} \circ \Phi[\etaup_p(\text{dist}(p,f^N(x)))] \left(A^N(x)\circ \pi_{E^u_{x,A}} \circ \Phi[\etaup_p(\text{dist}(p,x))] \cdot v^u\right)\|$$
$$\leq C_1C\beta^N \|A(f^N(x)) \circ \pi_{E^u_{f^N(x),A}} \circ \Phi[\etaup_p(\text{dist}(p,f^N(x)))]\left(A^N(x)\circ \pi_{E^u_{x,A}} \circ \Phi[\etaup_p(\text{dist}(p,x))] \cdot v^u\right)\|.$$

\bigskip

Notice that this contribution to the estimation of the unstable Lyapunov exponent corresponds to just one iterate of $f$ since, by the choice of $r$, now the orbit of $x$ leaves $B(p,r)$.

\bigskip

As the orbit of $x$ returns to $B(p,r)$ infinitely often but spends the longest periods out of this ball, we conclude that the unstable Lyapunov exponent strictly decreases by the amount $\ln(\Delta)$, and so the sum of center exponents must increase by the same quantity. Therefore,
$$
\sum n_{c,B} \lambda^c_{B} \simeq \sum n_{c,A} \lambda^c_{A}-\ln(\Delta) =-\ln(\Delta)>0.
$$
\end{subsubsection}

\bigskip

\begin{subsubsection}{\textbf{Second case: $d>1$}}

There are only a few adjustments to be done on the previous argument, the ones regarding the steps where we used that the dimension of $E^u$ is $1$. Namely:

\bigskip

\noindent (a) The evaluation of the sum of the unstable Lyapunov exponents.

\bigskip

\noindent According to (b.2) of Theorem \ref{Ruelle}, we may start with a unit vector $v^u \in V_i \setminus V_{i+1}$ and check that the corresponding Lyapunov exponent decreases. The reasoning proceeds precisely as before till the estimation (\ref{III}).

\bigskip

\noindent (b) The comparison of the two vectors in (\ref{III}).

\bigskip

\noindent The quotient
$$\frac{\|A(f^N(x))\circ \pi_{E^u_{f^N(x),A}} \circ \Phi[\etaup_p(\text{dist}(p,f^N(x)))] \cdot \left(A^N(x)\circ \pi_{E^c_{x,A}} \circ \Phi[\etaup_p(\text{dist}(p,x))] \cdot v^u\right)\|}{\|A(f^N(x))\circ \pi_{E^u_{f^N(x),A}} \circ \Phi[\etaup_p(\text{dist}(p,f^N(x)))] \cdot \left(A^N(x)\circ \pi_{E^u_{x,A}} \circ \Phi[\etaup_p(\text{dist}(p,x))] \cdot v^u\right)\|}$$
is harder to deal now with since the vectors
$$\vartheta_1=\pi_{E^u_{f^N(x),A}} \circ \Phi[\etaup_p(\text{dist}(p,f^N(x)))] \left(A^N(x)\circ \pi_{E^c_{x,A}} \circ \Phi[\etaup_p(\text{dist}(p,x))] \cdot v^u\right)$$
and
$$\vartheta_2=\pi_{E^u_{f^N(x),A}} \circ \Phi[\etaup_p(\text{dist}(p,f^N(x)))]\left(A^N(x)\circ \pi_{E^u_{x,A}} \circ \Phi[\etaup_p(\text{dist}(p,x))] \cdot v^u\right),$$
both live in $E^u_{f^N(x),A}$ but this space is no longer one-dimensional. However, as this concerns only one iterate of $A$ -- because, by the choice of $r$, the orbit of $x$ will leave $B(p,r)$ at once --, it is enough to check that, applying $A(f^N(x))$ to the inequality
$$\|\pi_{E^u_{f^N(x),A}} \circ \Phi[\etaup_p(\text{dist}(p,f^N(x)))] \left(A^N(x)\circ \pi_{E^c_{x,A}} \circ \Phi[\etaup_p(\text{dist}(p,x))] \cdot v^u\right)\| \leq$$ $$C_1C\beta^N \|\pi_{E^u_{f^N(x),A}} \circ \Phi[\etaup_p(\text{dist}(p,f^N(x)))]\left(A^N(x)\circ \pi_{E^u_{x,A}} \circ \Phi[\etaup_p(\text{dist}(p,x))] \cdot v^u\right)\|,$$
we get
$$\|A(f^N(x)) (\vartheta_1) \leq \|A\| \, \|\vartheta_1\| \leq \|A\| \, C_1C\beta^N \|\vartheta_2\|,$$
which ends the proof.
\end{subsubsection}

\end{proof}

\bigskip

If $1-\Delta$ is small enough, we may ensure that there is still an uniform gap between the Lyapunov exponents that correspond to different bundles of the Oseledets-Ruelle's decomposition of $B$. Therefore, as for $B$ the decomposition $E_B^{cu}=E_B^u \bigoplus E_B^s$ is finest, we deduce that
\begin{corollary}
$B$ is partially hyperbolic and $E_B^u=E_A^u$, $E_B^c=E_A^c$ and $E_B^s= E_A^s$.
\end{corollary}

\begin{proof}
If $\lambda^A_1 > \cdots > \lambda^A_d$ are the Lyapunov exponents associated with $E^u_A$ and $\lambda^A_{d+1} > \cdots \lambda^A_{D}$ the ones corresponding to $E^c_A$, it is enough to demand that $\lambda^A_d + \ln(\Delta) > \lambda^A_{d+1}$, that is, $\exp(\lambda^A_{d+1}-\lambda^A_d)< \Delta <1$, which amounts to consider a small rotation on the action of $\Phi \circ \etaup$ within Lemma \ref{Perturbation on the center-unstable space}.
\end{proof}
\end{subsection}

\medskip

\begin{subsection}{Perturbation on the central space}\label{PC}

Let $A\in C^{0}_{I}(X,\mathcal{C}(\mathcal{H}))$ be a partially hyperbolic cocycle with an extended Oseledets-Ruelle's decomposition $E_x^{u}\oplus E_x^{c}\oplus E_x^{s}=\mathcal{H}$, such that, for any $x \in \mathcal{O}(A)$, the space $E_{x,A}^{u}$ is non-trivial with dimension $d$ and the dimension of $E_x^{cu}$ is $D$.

For any $p \in \{1,2,\cdots,D\}$, consider the map
$$
\Lambda_p: C \in C^{0}_{I}(X,\mathcal{C}(\mathcal{H})) \mapsto \lambda^C_1+ \cdots + \lambda^C_p,
$$
where $\lambda^C_j$ is the $j$th Lyapunov exponent of the cocycle $C$. As $\Lambda_p$ is upper semicontinuous, for any $p$ (see section $3.5$ of \cite{BC}), defined on a Baire space (see section $3.1$ of \cite{BC}), it has a residual set $\mathcal{R}_p$ of continuity points. Therefore there is a partially hyperbolic cocycle $A_0$ inside the residual set
$$\mathcal{R}_1\cap \cdots \mathcal{R}_D \cap \mathcal{SF}_{E^u, E^c}$$
close enough to $A$ so that we are sure that its norm is positive. Thus, given $\epsilon >0$, there is a neighborhood $\mathcal{U}$ of $A_0$ such that (see Lemma \ref{Persistence}):
\begin{enumerate}
\item $\forall C \in \mathcal{U}$, $C$ is partially hyperbolic, $\text{dim}(E^u_C)= d$ and $\text{dim}(E^{cu}_C)= D.$
\item $\forall C_1, C_2 \in \mathcal{U} \, \, \forall p \in \{1,\cdots,D\} \,\, \left|\Lambda_p(C_1) - \Lambda_p(C_2)\right|<\epsilon.$
\end{enumerate}

\medskip

If $(\Lambda_D-\Lambda_d)(A_0) \neq 0$, take $B=A_0$ and proceed to the next paragraph. Otherwise, apply to $A_0$ the strategy presented on the previous subsection to get a cocycle $B \in \mathcal{U}$ whose sum of central Lyapunov exponents is positive.

If either $D=d$ or the Lyapunov exponents corresponding to $E_B^{c}$ are all equal, there is nothing else to be proved. If both conditions fail, take two distinct Lyapunov exponents, say $\lambda_p > \lambda_{p+1}$, in $E^c_B$. As the sum $E_B^{cu}=E_B^u \oplus E_B^c$ is dominated and finest, and both $\lambda_p$ and $\lambda_{p+1}$ belong to $E_B^{c}$, there is no dominated sum $E_B^{cu}=V_1 \oplus V_2$ with $\text{dim}(V_1)= p$. Therefore, as $\lambda_{p+1}\neq -\infty$, by Lemma 4.4 in \cite{BC}, there is a cocycle $C \in \mathcal{U}$ close to $B$ such that
$$
\Lambda_p(C)<\Lambda_p(B)-\frac{\lambda_p(B)-\lambda_{p+1}(B)}{2} + \epsilon.
$$
Hence
$$
\lambda_p(B)-\lambda_{p+1}(B)<2\left|\Lambda_p(C)-\Lambda_p(B)\right| + 2\epsilon<4\epsilon,
$$
which means that the central Lyapunov exponents of $B$ are all close to each other. Hence each one is approximately equal to $\frac{\Lambda_D - \Lambda_d}{D-d}$, and so does not vanish.
\end{subsection}
\end{section}

\bigskip

\section*{Acknowledgements}
MB was partially supported by FCT (Funda\c{c}\~ao para a Ci\^encia e a Tecnologia) through CMUP (SFRH/BPD/20890/2004) and the project PTDC/MAT/099493/2008.

\end{document}